\documentclass[a4paper,twoside,11pt]{article}
\usepackage{bbm}
\usepackage{amsfonts}
\usepackage{amssymb}\usepackage{amsmath}
 \numberwithin{equation}{section}
\begin{document}

\title{Blow-up phenomena and global existence for a periodic two-component Hunter-Saxton system}

\author{Jingjing Liu\footnote{e-mail:
jingjing830306@163.com }\\
Department of Mathematics,
Sun Yat-sen University,\\
510275 Guangzhou, China
\bigskip\\
Zhaoyang Yin\footnote{e-mail: mcsyzy@mail.sysu.edu.cn}
\\ Department of Mathematics, Sun Yat-sen University,\\ 510275 Guangzhou, China}
\date{}
\maketitle

\begin{abstract}
This paper is concerned with blow-up phenomena and global existence
for a periodic two-component Hunter-Saxton system. We first derive
the precise blow-up scenario for strong solutions to the system.
Then, we present several new blow-up results of strong solutions and
a new global existence result to the system. Our obtained results
for the system are sharp and improve considerably
earlier results.\\

\end{abstract}

\section{Introduction}
\par
In this paper, we study the Cauchy problem of the following periodic
two-component  Hunter-Saxton system:
\begin{equation}
 \ \ \ \ \ \ \ \  \ \ \ \ \left\{\begin{array}{ll}
u_{txx}+2u_{x}u_{xx}+uu_{xxx}-k\rho\rho_{x}=0 ,&t > 0,\,x\in \mathbb{R},\\
\rho_{t}+(\rho u)_x=0, &t > 0,\,x\in \mathbb{R},\\
u(0,x) = u_{0}(x),& x\in \mathbb{R}, \\
\rho(0,x) = \rho_{0}(x),&x\in \mathbb{R},\\
u(t,x+1)=u(t,x), & t \geq 0, x\in \mathbb{R},\\
\rho(t,x+1)=\rho(t,x), & t \geq 0, x\in \mathbb{R},\\ \end{array}\right. \\
\end{equation}
where $k=\pm 1$. The system (1.1) was originally proposed in
\cite{P-P} and is the short-wave limit of the two-component
Camassa-Holm system \cite{C-I, E-L-Y}. The system (1.1) is also a
special case of Green-Naghdi system modeling the non-dissipative
dark matter \cite{M-V-P}.

For $\rho\equiv 0$, the system (1.1) reduces to the Hunter-Saxton
equation \cite{J-R}, which describes the propagation of weakly
nonlinear orientation waves in a massive nematic liquid crystal
director field. The single-component model also arises in a
different physical context as the high-frequency limit \cite{HHD,
J-Y} of the Camassa-Holm equation for shallow water waves \cite{R-D,
J-H}, a re-expression of the geodesic flow on the diffeomorphism
group of the circle \cite{C-K} with a bi-Hamiltonian structure
\cite{F-F} which is completely integrable \cite{C-M}. The
Hunter-Saxton equation also has a bi-Hamiltonian structure
\cite{J-H, P-P} and is completely integrable \cite{B1, J-Y}.
Moreover, the Hunter-Saxton equation has a geometric interpretation
which was intensively studied in \cite{JL}.

The initial value problem for the Hunter-Saxton equation on the line
(nonperiodic case) was studied by Hunter and Saxton in \cite{J-R}.
Using the method of characteristics, they showed that smooth
solutions exist locally and break down in finite time, see
\cite{J-R}. The occurrence of blow-up can be interpreted physically
as the phenomenon by which waves that propagate away from the
perturbation ¡°knock¡± the director field out of its unperturbed
state \cite{J-R}. The initial value problem for the Hunter-Saxton
equation on the unit circle $\mathbb{S}=\mathbb{R}/\mathbb{Z}$ was
discussed in \cite{Y}. The author proved the local existence of
strong solutions to the periodic Hunter-Saxton equation, showed that
all strong solutions except space-independent solutions blow up in
finite time by using Kato semigroup method \cite{Kato 1}. Moreover,
the behavior of the solutions exhibits different features.

For $\rho\not\equiv0$, peakon solutions and the Cauchy problem of
the system (1.1) with $k=\pm 1$ have been discussed in \cite{C-I,
MW}. Recently, a generalization of the two-component Hunter-Saxton
system was proposed in \cite{MW2}. The global existence of solutions
to the generalized two-component Hunter-Saxton system was obtained
in \cite{hm}. The aim of this paper is to study further blow-up
phenomena and global existence of the system (1.1). The precise
blow-up scenario, several new blow-up results and a new global
existence result of strong solutions to the system (1.1) are
presented. The obtained results are sharp and improve considerably
the recent results in \cite{C-I, MW}.

The paper is organized as follows. In Section 2, we recall the local
existence of the initial value problem associated with the system
(1.1). In Section 3, we derive two precise blow-up scenarios. In
Section 4, we present several explosion criteria of strong solutions
to the system (1.1) with rather general initial data. In Section 5,
we give a new global existence result of strong solutions to the
system (1.1).

 \textbf{Notation}  Given a Banach space $Z$, we denote its norm by
 $\|\cdot\|_{Z}$. Since all space of functions are over
 $\mathbb{S}$, for simplicity, we drop $\mathbb{S}$ in our notations
  if there is no ambiguity. We let $[A,B]$ denote
 the commutator of linear operator $A$ and $B$. For convenience, we
 let $(\cdot|\cdot)_{s\times r}$ and $(\cdot|\cdot)_{s}$ denote the
 inner products of $H^{s}\times H^{r}$, $s,r\in \mathbb{R}_{+}$ and
 $H^{s}$, $s\in \mathbb{R}_{+}$, respectively.

\section{Local existence}
\newtheorem {remark2}{Remark}[section]
\newtheorem{theorem2}{Theorem}[section]
\newtheorem{lemma2}{Lemma}[section]

We provide now the framework in which we shall reformulate the
system (1.1). Integrating the first equation in (1.1) with respect
to $x$, we have
$$u_{tx}+uu_{xx}+\frac{1}{2}u_{x}^{2}-\frac{k}{2}\rho^{2}=a(t),$$
where
$$a(t)=-\frac{1}{2}\int_{\mathbb{S}}(k\rho^{2}+u_{x}^{2})dx$$ and $$\frac{d}{dt}a(t)=0,$$ cf. \cite{MW}. For convenience,
we let $a:=a(0)$. Thus,
\begin{equation}
u_{tx}+uu_{xx}=\frac{k}{2}\rho^{2}-\frac{1}{2}u_{x}^{2}+a.
\end{equation}
Integrating (2.1) with respect to $x$, we get
\begin{equation}
u_{t}+uu_{x}=\partial_{x}^{-1}(\frac{k}{2}\rho^{2}+\frac{1}{2}u_{x}^{2}+a)+h(t),
\end{equation}
where $\partial_{x}^{-1}g(x)=\int_{0}^{x}g(y)dy$ and $h(t):
[0,\infty)\rightarrow \mathbb{R}$ is an arbitrary continuous
function.

Thus we get an equivalent form of the system (1.1)
\begin{equation}
 \ \ \ \ \ \ \ \  \ \ \ \ \left\{\begin{array}{ll}
u_{t}+uu_{x}=\partial_{x}^{-1}(\frac{k}{2}\rho^{2}+\frac{1}{2}u_{x}^{2}+a)+h(t),&t > 0,\,x\in \mathbb{R},\\
\rho_{t}+(\rho u)_x=0, &t > 0,\,x\in \mathbb{R},\\
u(0,x) = u_{0}(x),&x\in \mathbb{R}, \\
\rho(0,x) = \rho_{0}(x),&x\in \mathbb{R},\\
\rho(t,x+1)=\rho(t,x), & t \geq 0, x\in \mathbb{R},\\
u(t,x+1)=u(t,x), & t \geq 0, x\in \mathbb{R},\\
\end{array}\right. \\
\end{equation}
where $k=\pm 1$, $\partial_{x}^{-1}g(x)=\int_{0}^{x}g(y)dy$ and
$h(t): [0,\infty)\rightarrow \mathbb{R}$ is an arbitrary continuous
function.

We now recall the local well-posedness result for system (2.3).

\begin{theorem2}{\cite{MW}}
Given $h(t)\in C([0,\infty); \mathbb{R})$ and
$z_{0}=(u_{0},\rho_{0})\in H^{s}\times H^{s-1}$, $s\geq 2,$ then
there exists a maximal $T=T(a, h(t), \parallel
z_{0}\parallel_{H^{s}\times H^{s-1}})>0$, and a unique solution
$z=(u, \rho)$ to (2.3) such that
$$ z=z(\cdot,z_{0})\in C([0,T); H^{s}\times H^{s-1})\cap
C^{1}([0,T);H^{s-1}\times H^{s-2}).
$$ Moreover, the solution depends continuously on the initial data, i.e., the
mapping $$z_{0}\rightarrow z(\cdot,z_{0}): H^{s}\times
H^{s-1}\rightarrow C([0,T); H^{s}\times H^{s-1})\cap
C^{1}([0,T);H^{s-1}\times H^{s-2})
$$
is continuous.
\end{theorem2}

As a consequence of Theorem 2.1 and the relation between the
solution of the system (1.1) and the solution of the system (2.3),
we have the following local exist result.
\begin{theorem2}
Given $z_{0}=(u_{0},\rho_{0})\in H^{s}\times H^{s-1}, s\geq 2.$ Then
there exists locally a family of solutions to (1.1).
\end{theorem2}
Note that the solution of the system (2.3) for any fixed $h(t)$ is
unique. However, the solution of the system (1.1) given by Theorem
2.2 is not unique by the arbitrariness of $h(t)$. In the following
sections, we discuss the corresponding unique solution to the system
(2.3) with a fixed $h(t)$.
\section{The precise blow-up scenario}
\newtheorem {remark3}{Remark}[section]
\newtheorem{theorem3}{Theorem}[section]
\newtheorem{lemma3}{Lemma}[section]
\newtheorem{definition3}{Definition}[section]
\newtheorem{claim3}{Claim}[section]

In this section, we present the precise blow-up scenario for strong
solutions to the system (1.1).

We first recall the following lemmas.

\begin{lemma3}
\cite{Kato 3} If $r>0$, then $H^{r}\cap L^{\infty}$ is an algebra.
Moreover
$$\parallel fg\parallel_{H^r}\leq c(\parallel  f\parallel_{L^{\infty}}\parallel g\parallel_{H^r}+\parallel
  f\parallel_{H^r}\parallel g\parallel_{L^{\infty}}),
$$
where c is a constant depending only on r.
\end{lemma3}

\begin{lemma3}
\cite{Kato 3} If $r>0$, then
$$\parallel[\Lambda^{r},f]g\parallel_{L^{2}}\leq c(\parallel \partial_{x}f\parallel_{L^{\infty}}
\parallel\Lambda^{r-1}g\parallel_{L^2}+\parallel
\Lambda^r f\parallel_{L^2}\parallel g\parallel_{L^{\infty}}),
$$
where c is a constant depending only on r.
\end{lemma3}

\begin{lemma3}\cite{Constantin 4}
Let $t_{0}>0$ and $v\in C^{1}([0,t_{0}); H^{2}(\mathbb{R}))$. Then
for every $t\in[0,t_{0})$ there exists at least one point $\xi(t)\in
\mathbb{R}$ with
$$ m(t):=\inf_{x\in \mathbb{R}}\{v_{x}(t,x)\}=v_{x}(t,\xi(t)),$$ and
the function $m$ is almost everywhere differentiable on $(0,t_{0})$
with $$ \frac{d}{dt}m(t)=v_{tx}(t,\xi(t)) \ \ \ \ a.e.\ on \
(0,t_{0}).$$
\end{lemma3}

\begin{remark3}
If $v\in C^{1}([0,t_{0}); H^{s}(\mathbb{R})), \ s>\frac{3}{2},$ then
Lemma 3.3 also holds true. Meanwhile, Lemma 3.3 works analogously
for $$M(t):= \sup\limits_{x\in\mathbb{R}}\{{v_{x}(t,x)}\}.$$
\end{remark3}

\begin{lemma3}
Assume $k=1$. Let $z_{0}=\left(\begin{array}{c}
                                u_{0} \\
                                \rho_{0} \\
                              \end{array}
                            \right)
\in H^{s}\times H^{s-1}$, $s > \frac{5}{2}$, be given and assume
that T is the maximal existence time of the corresponding solution
                  $z=\left(\begin{array}{c}
                                      u \\
                                      \rho \\
                                    \end{array}
                                  \right)$
to (1.1) with the initial data $z_{0}$. Then
$$\|\rho_{x}(t,\cdot)\|_{L^{\infty}}\leq K \exp\left\{-2\int_{0}^{t}u_{x}(s,\xi(s))ds\right\},$$ where
$(s,\xi(s))$ is a maximal point of $u_{xx}^{2}+\rho_{x}^{2}$ in
$[0,T)\times\mathbb{S}$ and
$K=\|u_{0,xx}\|_{L^{\infty}}+\|\rho_{0,x}\|_{L^{\infty}}.$
\end{lemma3}
\textbf{Proof} Multiplying the first equation in (1.1) by $u_{xx}$,
we get
\begin{equation}
\frac{1}{2}(u_{xx}^{2})_{t}+2u_{x}u_{xx}^{2}+u\frac{1}{2}(u_{xx}^{2})_{x}-\rho\rho_{x}u_{xx}=0.
\end{equation}
Differentiating the second equation in (1.1) in $x$ and multiplying
the obtained equation by $\rho_{x}$, we get
\begin{equation}
\frac{1}{2}(\rho_{x}^{2})_{t}+2u_{x}\rho_{x}^{2}+u\frac{1}{2}(\rho_{x}^{2})_{x}+\rho\rho_{x}u_{xx}=0.
\end{equation}
Adding the above two equations, we have
\begin{equation}
\frac{1}{2}(u_{xx}^{2}+\rho_{x}^{2})_{t}+2u_{x}(u_{xx}^{2}+\rho_{x}^{2})+\frac{1}{2}(u_{xx}^{2}+\rho_{x}^{2})_{x}=0.
\end{equation}
By $z\in C([0,T);H^{s}\times H^{s-1}),$ $s> \frac{5}{2}$, we know
$u_{xx}\in H^{s-2}$, $\rho_{x}\in H^{s-2}$. Moreover, since
$H^{s-2}$ is a Banach algebra for $s>\frac{5}{2},$
$u_{xx}^{2}+\rho_{x}^{2}\in H^{s-2},$ $s>\frac{5}{2}.$ Let
$M(t)=\sup\limits_{x\in\mathbb{S}}(u_{xx}^{2}+\rho_{x}^{2})(t,x).$
It follows from Remark 3.1 that there is a point $(t,\xi(t))\in
[0,T)\times \mathbb{S}$ such that
$M(t)=(u_{xx}^{2}+\rho_{x}^{2})(t,\xi(t)).$ Evaluating (3.3) on $(t,
\xi(t))$ we get
$$\frac{dM(t)}{dt}=-4u_{x}(t,\xi(t))M(t).$$ Then, we obtain
$$M(t)=M(0)\exp\left\{\int_{0}^{t}-4u_{x}(s,\xi(s))ds\right\}.$$
Note that
$$M(0)=\sup\limits_{x\in\mathbb{S}}(u_{0,xx}^{2}+\rho_{0,x}^{2})\leq
\|u_{0,xx}\|_{L^{\infty}}^{2}+\|\rho_{0,x}\|_{L^{\infty}}^{2}.$$
Thus, we get
$$\|\rho_{x}\|_{L^{\infty}}\leq K\exp\left\{-2\int_{0}^{t}u_{x}(s,\xi(s))ds\right\}.$$

Next we prove the following useful result on global existence of
solutions to (1.1).
\begin{theorem3} Assume $k=1$. Let $z_{0}=\left(\begin{array}{c}
                                u_{0} \\
                                \rho_{0} \\
                              \end{array}
                            \right)
\in H^{s}\times H^{s-1}$, $s > \frac{5}{2}$, be given and assume
that T is the maximal existence time of the corresponding solution
                  $z=\left(\begin{array}{c}
                                      u \\
                                      \rho \\
                                    \end{array}
                                  \right)$
to (2.3) with the initial data $z_{0}$. If there exists $M>0$ such
that
$$
\| u_{x}(t,\cdot)\|_{L^{\infty}}
+\|\rho(t,\cdot)\|_{_{L^{\infty}}}\leq M,\ \ t\in[0,T),
$$
then the $H^s\times H^{s-1}$-norm of $z(t,\cdot)$ does not blow up
on [0,T).
 \end{theorem3}
\textbf{Proof} Let $z=\left(
                                    \begin{array}{c}
                                      u \\
                                      \rho \\
                                    \end{array}
                                  \right)$ be the solution to
(2.3) with the initial data $z_{0}\in H^s\times H^{s-1},\
s>\frac{5}{2}$, and let T be the maximal existence time of the
corresponding solution $z$, which is guaranteed by Theorem 2.1.
Throughout this proof, $c>0$ stands for a generic constant depending
only on $s$.

By $\|u_{x}(t,\cdot)\|_{L^{\infty}}\leq M$ and Lemma 3.4, we get
\begin{equation}
\|\rho_{x}\|_{L^{\infty}}
\leq K\exp\left\{2\int_{0}^{t}u_{x}(s,\xi(s))ds\right\}\\
\leq Ke^{2Mt}:=c(t).
\end{equation}

Applying the operator $\Lambda^{s}$ to the first equation in (2.3),
multiplying by $\Lambda^{s} u$, and integrating over $\mathbb{S}$,
we obtain
\begin{equation}
\frac{d}{dt}\|u\|^{2}_{H^{s}}=-2(uu_{x},
u)_{s}+2(u,\partial_{x}^{-1}(\frac{1}{2}\rho^{2}+\frac{1}{2}u_{x}^{2}+a)+h(t))_{s}.
\end{equation}
Let us estimate the first term of the right-hand side of (3.5).
\begin{align*}
|(uu_{x},u)_{s}|&= |(\Lambda^s(u\partial_{x}u),\Lambda^s u)_{0}|\\
&=|([\Lambda^s,u]\partial_{x}u,\Lambda^s
u)_{0}+(u\Lambda^s\partial_{x}u,\Lambda^s u)_{0}|\\
&\leq \|[\Lambda^s,u]\partial_{x}u\|_{L^2}\|\Lambda^s
u\|_{L^2}+\frac{1}{2}|(u_{x}\Lambda^s u,\Lambda^s u)_{0}|\\
& \leq (c\| u_{x}\|_{L^{\infty}}+\frac{1}{2}\|
u_{x}\|_{L^{\infty}})\|u\|^2_{H^s}\\
&\leq c\|u_{x}\|_{L^{\infty}}\|u\|^2_{H^s},
\end{align*}
where we used Lemma 3.2 with $r=s$. Then, we estimate the second
term of the right-hand side of (3.5) in the following way:
\begin{align*}
&|(\partial_{x}^{-1}(\frac{1}{2}\rho^{2}+\frac{1}{2}u_{x}^{2}+a)+h(t),u)_{s}|\\
\leq \
&\|\partial_{x}^{-1}(\frac{1}{2}\rho^{2}+\frac{1}{2}u_{x}^{2}+a)+h(t)\|_{H^{s}}\|u\|_{H^{s}}\\
\leq  \
&(\|\partial_{x}^{-1}(\frac{1}{2}\rho^{2}+\frac{1}{2}u_{x}^{2}+a)\|_{L^{2}}+\|\frac{1}{2}\rho^{2}+\frac{1}{2}u_{x}^{2}+a\|_{H^{s-1}}
+\|h(t)\|_{H^{s}})\|u\|_{H^{s}}\\
\leq \
&(\|\frac{1}{2}\rho^{2}+\frac{1}{2}u_{x}^{2}+a\|_{L^{2}}+\|\frac{1}{2}\rho^{2}+\frac{1}{2}u_{x}^{2}+a\|_{H^{s-1}}
+\|h(t)\|_{H^{s}})\|u\|_{H^{s}}\\
\leq  \
&c(\|\rho^{2}\|_{H^{s-1}}+\|u_{x}^{2}\|_{H^{s-1}}+2\|a\|_{H^{s-1}}+\|h(t)\|_{H^{s}})\|u\|_{H^{s}}\\
\leq  \ &c
(\|\rho\|_{L^{\infty}}\|\rho\|_{H^{s-1}}+\|u_{x}\|_{L^{\infty}}\|u_{x}\|_{H^{s-1}}+|a|+|h(t)|)\|u\|_{H^{s}}\\
\leq \  &c
(\|\rho\|_{L^{\infty}}+\|u_{x}\|_{L^{\infty}}+1)(\|u\|_{H^{s}}^{2}+\|\rho\|_{H^{s-1}}^{2}+1),
\end{align*}
where we used Lemma 3.1 with $r=s-1$. Combining the above two
inequalities with (3.5), we get
\begin{equation}
\frac{d}{dt}\|u\|^{2}_{H^{s}}\leq
c(\|\rho\|_{L^{\infty}}+\|u_{x}\|_{L^{\infty}}+1)(\|u\|_{H^{s}}^{2}+\|\rho\|_{H^{s-1}}^{2}+1).
\end{equation}
In order to derive a similar estimate for the second component
$\rho$, we apply the operator $\Lambda^{s-1}$ to the second equation
in (2.3), multiply by $\Lambda^{s-1}\rho$, and integrate over
$\mathbb{S}$, to obtain
\begin{equation}
\frac{d}{dt}\|\rho\|_{H^{s-1}}^{2}=-2(u\rho_{x},\rho)_{s-1}-2(u_{x}\rho,\rho)_{s-1}.
\end{equation}
Let us estimate the first term of the right hand side of (3.7)
\begin{align*}
&|(u\rho_{x},\rho)_{s-1}|\\
=  \ &|(\Lambda^{s-1}(u\partial_{x}\rho),\Lambda^{s-1}\rho)_{0}|\\
=  \ &|([\Lambda^{s-1},u]\partial_{x}\rho,
\Lambda^{s-1}\rho)_{0}+(u\Lambda^{s-1}\partial_{x}\rho,
\Lambda^{s-1}\rho)_{0}|\\
\leq \
&\|[\Lambda^{s-1},u]\partial_{x}\rho\|_{L^{2}}\|\Lambda^{s-1}\rho\|_{L^{2}}+\frac{1}{2}|(u_{x}\Lambda^{s-1}\rho,\Lambda^{s-1}\rho)_{0}|\\
\leq \
&c(\|u_{x}\|_{L^{\infty}}\|\rho\|_{H^{s-1}}+\|\rho_{x}\|_{L^{\infty}}\|u\|_{H^{s-1}})\|\rho\|_{H^{s-1}}
+\frac{1}{2}\|u_{x}\|_{L^{\infty}}\|\rho\|_{H^{s-1}}^{2}\\
\leq \
&c(\|u_{x}\|_{L^{\infty}}+\|\rho_{x}\|_{L^{\infty}})(\|\rho\|_{H^{s-1}}^{2}+\|u\|_{H^{s}}^{2}),
\end{align*}
here we applied Lemma 3.2 with $r=s-1$. Then we estimate the second
term of the right hand side of (3.7). Based on Lemma 3.1 with
$r=s-1$, we get
\begin{align*}
|(u_{x}\rho,\rho)_{s-1}|\leq \ &\|u_{x}\rho\|_{H^{s-1}}\|\rho\|_{H^{s-1}}\\
\leq \
&c(\|u_{x}\|_{L^{\infty}}\|\rho\|_{H^{s-1}}+\|\rho\|_{L^{\infty}}\|u_{x}\|_{H^{s-1}})\|\rho\|_{H^{s-1}}\\
\leq \
&c(\|u_{x}\|_{L^{\infty}}+\|\rho_{x}\|_{L^{\infty}})(\|\rho\|_{H^{s-1}}^{2}+\|u\|_{H^{s}}^{2}).
\end{align*}
Combining the above two inequalities with (3.7), we get
\begin{equation}
\frac{d}{dt}\|\rho\|^{2}_{H^{s-1}}\leq
c(\|u_{x}\|_{L^{\infty}}+\|\rho\|_{L^{\infty}}+\|\rho_{x}\|_{L^{\infty}})(\|u\|_{H^{s}}^{2}+\|\rho\|_{H^{s-1}}^{2}+1).
\end{equation}
By (3.6) and (3.8), we have
\begin{align*}
&\frac{d}{dt}(\|\rho\|^{2}_{H^{s-1}}+\|u\|_{H^{s}}^{2}+1)\\
\leq \
&c(\|u_{x}\|_{L^{\infty}}+\|\rho\|_{L^{\infty}}+\|\rho_{x}\|_{L^{\infty}}+1)(\|u\|_{H^{s}}^{2}+\|\rho\|_{H^{s-1}}^{2}+1).
\end{align*}
An application of (3.4), Gronwall's inequality and the assumption of
the theorem yield
$$(\|\rho\|^{2}_{H^{s-1}}+\|u\|_{H^{s}}^{2}+1)\leq\exp(c(M+c(t)+1))(\|\rho_{0}\|^{2}_{H^{s-1}}+\|u_{0}\|_{H^{s}}^{2}+1).$$
This completes the proof of the theorem.

Given $z_{0}\in H^{s}\times H^{s-1}$ with $s\geq 2$. Theorem 2.1
ensures the existence of a maximal  $T> 0$ and a solution
                  $z=\left(\begin{array}{c}
                                      u \\
                                      \rho \\
                                    \end{array}
                                  \right)$
to (2.3) such that
$$
z=z(\cdot,z_{0})\in C([0,T); H^{s}\times H^{s-1})\cap
C^{1}([0,T);H^{s-1}\times H^{s-2}).
$$

Consider now the following initial value problem

\begin{equation}
\left\{\begin{array}{ll}q_{t}=u(t,q),\ \ \ \ t\in[0,T), \\
q(0,x)=x,\ \ \ \ x\in\mathbb{R}, \end{array}\right.
\end{equation}
where $u$ denotes the first component of the solution $z$ to (2.3).
Then we have the following two useful lemmas.

Applying classical results in the theory of ordinary differential
equations, one can obtain the following result on $q$ which is
crucial in the proof of blow-up scenarios.

\begin{lemma3}\cite{E-L-Y, gm}
 Let $u\in C([0,T); H^{s})\bigcap C^{1}([0,T); H^{s-1}), s\geq 2$. Then Eq.(3.9) has a unique solution
$q\in C^1([0,T)\times \mathbb{R};\mathbb{R})$. Moreover, the map
$q(t,\cdot)$ is an increasing diffeomorphism of $\mathbb{R}$ with
$$
q_{x}(t,x)=exp\left(\int_{0}^{t}u_{x}(s,q(s,x))ds\right)>0, \ \
(t,x)\in [0,T)\times \mathbb{R}.$$
\end{lemma3}
Following the similar proof in \cite{E-L-Y}, we obtain the next
result:
\begin{lemma3}
Assume $k= \pm 1$. Let $z_{0}=\left(\begin{array}{c}
                                u_{0} \\
                                \rho_{0} \\
                              \end{array}
                            \right)
\in H^{s}\times H^{s-1}$, $s \geq 2$ and let $T>0$ be the maximal
existence time of the corresponding solution
                  $z=\left(\begin{array}{c}
                                      u \\
                                      \rho \\
                                    \end{array}
                                  \right)$
to (1.1). Then we have
\begin{equation}
\rho(t,q(t,x))q_{x}(t,x)=\rho_{0}(x), \ \ \ \forall \
(t,x)\in[0,T)\times \mathbb{S}.
\end{equation}
Moreover, if there exists $M_{1}> 0$ such that $u_{x}\geq -M_{1}$
for all $(t,x)\in [0,T)\times \mathbb{S}$, then
$$\|\rho(t,\cdot)\|_{L^{\infty}}=\|\rho(t,q(t,\cdot))\|_{L^{\infty}}\leq
e^{M_{1}T}\|\rho_{0}(\cdot)\|_{L^{\infty}}, \ \ \ \forall \ t\in[0,
T).$$ Furthermore, if $\rho_{0}\in L^{1}$, then
$$\int_{\mathbb{S}}|\rho(t,x)|dx=\int_{\mathbb{S}}|\rho_{0}(x)|dx, \
\ \forall \ t\in [0,T).$$
\end{lemma3}
Our next result describes the precise blow-up scenarios for
sufficiently regular solutions to (1.1).

\begin{theorem3}
Assume $k=1$. Let $z_{0}=\left(\begin{array}{c}
                                u_{0} \\
                                \rho_{0} \\
                              \end{array}
                            \right)
\in H^{s}\times H^{s-1}$, $s>\frac{5}{2} $ be given and let T be the
maximal existence time of the corresponding solution
                  $z=\left(\begin{array}{c}
                                      u \\
                                      \rho \\
                                    \end{array}
                                  \right)
$ to (2.3) with the initial data $z_{0}$. Then the corresponding
solution blows up in finite time if and only if
$$\liminf\limits_{t\rightarrow T}\inf\limits_{x\in \mathbb{S}}{u_{x}(t,x)}=-\infty. $$
\end{theorem3}

\textbf{Proof} By Theorem 2.1 and Sobolev's imbedding theorem it is
clear that if
$$\liminf\limits_{t\rightarrow T}\inf\limits_{x\in \mathbb{S}}{u_{x}(t,x)}=-\infty,$$ then $T<\infty$.

Let $T<\infty$. Assume that there exists $M_{1}>0$ such that
$$u_{x}(t,x)\geq -M_{1}, \ \ \forall\ (t,x)\in
[0,T)\times\mathbb{S}.$$ By Lemma 3.6, we
have$$\|\rho(t,\cdot)\|_{L^{\infty}}\leq
e^{M_{1}T}\|\rho_{0}\|_{L^{\infty}}, \ \ \forall\ t\in [0,T).$$

Take $K_{0}=K^{2}e^{4M_{1}T}.$ By the first equation in (2.3), a
direct computation implies the following inequality
\begin{align}
&\frac{d}{dt}\int_{\mathbb{S}}u(t,x)^{2}dx\\
\nonumber=\
&2\int_{\mathbb{S}}u\partial_{x}^{-1}(\frac{1}{2}u_{x}^{2}+\frac{k}{2}\rho^{2}+a)dx+2h(t)\int_{\mathbb{S}}u
dx\\
\nonumber\leq \
&\int_{\mathbb{S}}u^{2}dx+\frac{1}{4}\int_{\mathbb{S}}\left(\int_{0}^{x}(u_{y}^{2}+k\rho^{2}+2a)dy\right)^{2}dx+|h(t)|\left(1+\int_{\mathbb{S}}u(t,x)^{2}dx\right)\\
\nonumber\leq \
&|h(t)|+(1+|h(t)|)\int_{\mathbb{S}}u(t,x)^{2}dx+\frac{1}{4}\left(\int_{0}^{1}(u_{x}^{2}+\rho^{2}+2|a|)dx\right)^{2}\\
\nonumber= \
&|h(t)|+(1+|h(t)|)\int_{\mathbb{S}}u(t,x)^{2}dx+\frac{1}{4}\left[2|a|+\int_{0}^{1}(u_{0,x}^{2}+\rho_{0}^{2})dx\right]^{2}
\end{align}
for $t\in(0,T)$.

Multiplying (2.1) by $u_{x}$ and integrating by parts, we get
\begin{align}
\frac{d}{dt}\int_{\mathbb{S}}u_{x}^{2}dx= \ &2\int_{\mathbb{S}}u_{x}(-uu_{xx}+\frac{k}{2}\rho^{2}-\frac{1}{2}u_{x}^{2}+a)dx\\
\nonumber = \
&\int_{\mathbb{S}}-2uu_{x}u_{xx}dx+k\int_{\mathbb{S}}u_{x}\rho^{2}dx-\int_{\mathbb{S}}u_{x}^{3}dx+2a\int_{\mathbb{S}}u_{x}dx\\
\nonumber =\ &k\int_{\mathbb{S}}u_{x}\rho^{2}dx\\
\nonumber \leq \
&\|\rho\|_{L^{\infty}}^{2}+\|\rho\|_{L^{\infty}}^{2}\int_{\mathbb{S}}u_{x}^{2}dx.
\end{align}

Multiplying the first equation in (1.1) by $m=u_{xx}$ and
integrating by parts, we find
\begin{align}
\frac{d}{dt}\int_{\mathbb{S}}m^{2}dx= \ &-4\int_{\mathbb{S}}u_{x}m^{2}dx-2\int_{\mathbb{S}}umm_{x}dx+2k\int_{\mathbb{S}}m\rho\rho_{x}dx\\
\nonumber= \
&-3\int_{\mathbb{S}}u_{x}m^{2}dx+2k\int_{\mathbb{S}}m\rho\rho_{x}dx\\
\nonumber\leq \
&3M_{1}\int_{\mathbb{S}}m^{2}dx+\|\rho\|_{L^{\infty}}\int_{\mathbb{S}}m^{2}+\rho_{x}^{2}dx\\
\nonumber\leq \
&(3M_{1}+\|\rho\|_{L^{\infty}})\int_{\mathbb{S}}m^{2}dx+\|\rho\|_{L^{\infty}}\int_{\mathbb{S}}\rho_{x}^{2}dx.
\end{align}

Differentiating the first equation in (1.1) with respect to $x$,
multiplying the obtained equation by $m_{x}=u_{xxx},$ integrating by
parts and using Lemma 3.4, we obtain
\begin{align}
&\frac{d}{dt}\int_{\mathbb{S}}m_{x}^{2}dx\\
\nonumber= \
&-4\int_{\mathbb{S}}m^{2}m_{x}dx-6\int_{\mathbb{S}}u_{x}m_{x}^{2}-2\int_{\mathbb{S}}um_{xx}m_{x}
+2k\int_{\mathbb{S}}\rho_{x}^{2}m_{x}+2k\int_{\mathbb{S}}\rho\rho_{xx}m_{x}dx\\
\nonumber= \
&-5\int_{\mathbb{S}}u_{x}m_{x}^{2}dx+2k\int_{\mathbb{S}}\rho_{x}^{2}m_{x}dx+2k\int_{\mathbb{S}}\rho\rho_{xx}m_{x}dx\\
\nonumber\leq \
&5M_{1}\int_{\mathbb{S}}m_{x}^{2}dx+2\|\rho_{x}\|_{L^{\infty}}^{2}\int_{\mathbb{S}}|m_{x}|dx+\|\rho\|_{L^{\infty}}\int_{\mathbb{S}}(\rho_{xx}^{2}+m_{x}^{2})dx\\
\nonumber\leq \
&5M_{1}\int_{\mathbb{S}}m_{x}^{2}dx+\|\rho\|_{L^{\infty}}\int_{\mathbb{S}}(\rho_{xx}^{2}+m_{x}^{2})dx+2\|\rho_{x}\|_{L^{\infty}}^{2}
+2\|\rho_{x}\|_{L^{\infty}}^{2}\int_{\mathbb{S}}m_{x}^{2}dx\\
\nonumber\leq \
&(5M_{1}+\|\rho\|_{L^{\infty}}+2K_{0})\int_{\mathbb{S}}m_{x}^{2}dx+
\|\rho\|_{L^{\infty}}\int_{\mathbb{S}}\rho_{xx}^{2}dx+2K_{0}.
\end{align}

Multiplying the second equation in (1.1) by $\rho$ and integrating
by parts, we have
\begin{equation}
\frac{d}{dt}\int_{\mathbb{S}}\rho^{2}dx=-\int_{\mathbb{S}}u_{x}\rho^{2}dx\leq
M_{1}\int_{\mathbb{S}}\rho^{2}dx.
\end{equation}

Differentiating the second equation in (1.1) with respect to $x$,
multiplying the obtained equation by $\rho_{x}$ and integrating by
parts, we obtain
\begin{align}
\frac{d}{dt}\int_{\mathbb{S}}\rho_{x}^{2}dx= \
&-3\int_{\mathbb{S}}u_{x}\rho_{x}^{2}dx-2\int_{\mathbb{S}}m\rho\rho_{x}dx\\
\nonumber\leq \
&3M_{1}\int_{\mathbb{S}}\rho_{x}^{2}dx+\|\rho\|_{L^{\infty}}\int_{\mathbb{S}}(m^{2}+\rho_{x}^{2})dx\\
\nonumber\leq \ &
(3M_{1}+\|\rho\|_{L^{\infty}})\int_{\mathbb{S}}\rho_{x}^{2}dx+\|\rho\|_{L^{\infty}}\int_{\mathbb{S}}m^{2}dx.
\end{align}

Differentiating the second equation in (1.1) with respect to $x$
twice, multiplying the obtained equation by $\rho_{xx},$ integrating
by parts and using Lemma 3.4, we obtain
\begin{align}
&\frac{d}{dt}\int_{\mathbb{S}}\rho_{xx}^{2}dx\\
\nonumber= \
&-5\int_{\mathbb{S}}u_{x}\rho_{xx}^{2}dx+\int_{\mathbb{S}}u_{xxx}(3\rho_{x}^{2}-2\rho\rho_{xx})dx\\
\nonumber\leq \
&5M_{1}\int_{\mathbb{S}}\rho_{xx}^{2}dx+\int_{\mathbb{S}}m_{x}(3\rho_{x}^{2}-2\rho\rho_{xx})dx\\
\nonumber\leq \
&5M_{1}\int_{\mathbb{S}}\rho_{xx}^{2}dx+3\|\rho_{x}\|_{L^{\infty}}^{2}\int_{\mathbb{S}}|m_{x}|dx+
\|\rho\|_{L^{\infty}}\int_{\mathbb{S}}2m_{x}\rho_{xx}dx\\
\nonumber\leq \
&(5M_{1}+\|\rho\|_{L^{\infty}})\int_{\mathbb{S}}\rho_{xx}^{2}dx+(3\|\rho_{x}\|_{L^{\infty}}^{2}
+\|\rho\|_{L^{\infty}})\int_{\mathbb{S}}m_{x}^{2}+3\|\rho_{x}\|_{L^{\infty}}^{2}\\
\nonumber\leq \
&(5M_{1}+\|\rho\|_{L^{\infty}})\int_{\mathbb{S}}\rho_{xx}^{2}dx
+(3K_{0}+\|\rho\|_{L^{\infty}})\int_{\mathbb{S}}m_{x}^{2}dx+3K_{0}.
\end{align}

Summing (3.10)-(3.16), we have
\begin{align*}
&\frac{d}{dt}\int_{\mathbb{S}}(u^{2}+u_{x}^{2}+m^{2}+m_{x}^{2}+\rho^{2}+\rho_{x}^{2}+\rho_{xx}^{2})dx\\
\leq \
&K_{1}\int_{\mathbb{S}}(u^{2}+u_{x}^{2}+m^{2}+m_{x}^{2}+\rho^{2}+\rho_{x}^{2}+\rho_{xx}^{2})dx+K_{2},\\
\end{align*}
where
$$K_{1}=1+\max\limits_{t\in[0,T]}|h(t)|+8e^{M_{1}T}\|\rho_{0}\|_{L^{\infty}}+(e^{M_{1}T}\|\rho_{0}\|_{L^{\infty}})^{2}
+17M_{1}+5K_{0},$$
$$K_{2}=\max\limits_{t\in[0,T]}|h(t)|+\frac{1}{4}\left[2|a|+\int_{0}^{1}(u_{0,x}^{2}
+\rho_{0}^{2})dx\right]^{2}+(e^{M_{1}T}\|\rho_{0}\|_{L^{\infty}})^{2}+5K_{0}.$$
By means of Gronwall's inequality and the above inequality, we
deduce that
\begin{align*}
&\|u(t,\cdot)\|_{H^{3}}^{2}+\|\rho(t,\cdot)\|_{H^{2}}^{2}\\
\leq \
&e^{K_{1}t}(\|u_{0}\|_{H^{3}}^{2}+\|\rho_{0}\|_{H^{2}}^{2}+\frac{K_{2}}{K_{1}}),
\ \ \ \forall \ t\in[0,T).
\end{align*}
The above inequality, Sobolev's imbedding theorem and Theorem 3.1
ensure that the solution $z$ does not blow-up in finite time. This
completes the proof of the theorem.

Note that when $k=-1$, we cannot get Lemma 3.4. However, following
the similar proof of Theorems 3.1-3.2 we obtain the following two
results:
\begin{theorem3} Assume $k=-1$. Let $z_{0}=\left(\begin{array}{c}
                                u_{0} \\
                                \rho_{0} \\
                              \end{array}
                            \right)
\in H^{s}\times H^{s-1}$, $s \geq 2$, be given and assume that T is
the maximal existence time of the corresponding solution
                  $z=\left(\begin{array}{c}
                                      u \\
                                      \rho \\
                                    \end{array}
                                  \right)$
to (2.3) with the initial data $z_{0}$. If there exists $M>0$ such
that
$$
\| u_{x}(t,\cdot)\|_{L^{\infty}}
+\|\rho(t,\cdot)\|_{_{L^{\infty}}}+\|\rho_{x}(t,\cdot)\|_{L^{\infty}}\leq
M,\ \ t\in[0,T),
$$
then the $H^s\times H^{s-1}$-norm of $z(t,\cdot)$ does not blow up
on [0,T).
 \end{theorem3}

\begin{theorem3}
Assume $k=-1$. Let $z_{0}=\left(\begin{array}{c}
                                u_{0} \\
                                \rho_{0} \\
                              \end{array}
                            \right)
\in H^{s}\times H^{s-1}$, $s> \frac{5}{2}, $ be given and let T be
the maximal existence time of the corresponding solution
                  $z=\left(\begin{array}{c}
                                      u \\
                                      \rho \\
                                    \end{array}
                                  \right)
$ to (2.3) with the initial data $z_{0}$. Then the corresponding
solution blows up in finite time if and only if
$$\liminf\limits_{t\rightarrow T}\inf\limits_{x\in
\mathbb{S}}{u_{x}(t,x)}=-\infty \ \ \text{or} \ \
\limsup\limits_{t\rightarrow
T}\{\|\rho_{x}\|_{L^{\infty}}\}=+\infty.$$
\end{theorem3}

For initial data $z_{0}=\left(\begin{array}{c}
                                u_{0} \\
                                \rho_{0} \\
                              \end{array}
                            \right)
\in H^{2}\times H^{1}$, we have the following precise blow-up
scenario.

\begin{theorem3}
Assume $k=\pm 1$. Let $z_{0}=\left(\begin{array}{c}
                                u_{0} \\
                                \rho_{0} \\
                              \end{array}
                            \right)
\in H^{2}\times H^{1}$, and let T be the maximal existence time of
the corresponding solution
                  $z=\left(\begin{array}{c}
                                      u \\
                                      \rho \\
                                    \end{array}
                                  \right)
$ to (2.3) with the initial data $z_{0}$. Then the corresponding
solution blows up in finite time if and only if
$$\liminf\limits_{t\rightarrow T}\inf\limits_{x\in \mathbb{S}}{u_{x}(t,x)}=-\infty.$$
\end{theorem3}

\textbf{Proof} Let $z=\left(\begin{array}{c}
                                      u \\
                                      \rho \\
                                    \end{array}
                                  \right)
$ be the solution to (2.3) with the initial data $z_{0}\in
H^{2}\times H^{1},$ and let $T$ be the maximal existence time of the
solution $z$, which is guaranteed by Theorem 2.1.

Let $T< \infty$. Assume that there exists $M_{1}>0$ such that
$$u_{x}(t,x)\geq -M_{1}, \ \ \forall\ (t,x)\in
[0,T)\times\mathbb{S}.$$ By Lemma 3.6, we
have$$\|\rho(t,\cdot)\|_{L^{\infty}}\leq
e^{M_{1}T}\|\rho_{0}\|_{L^{\infty}}, \ \ \forall\ t\in [0,T).$$

Combining (3.11)-(3.13) and (3.15)-(3.16), we
obtain$$\frac{d}{dt}\int_{\mathbb{S}}u^{2}+u_{x}^{2}+m^{2}+\rho^{2}+\rho_{x}^{2}dx\leq
K_{3}\int_{\mathbb{S}}u^{2}+u_{x}^{2}+m^{2}+\rho^{2}+\rho_{x}^{2}dx+K_{4},$$
where
$$K_{3}=1+\max\limits_{t\in[0,T]}|h(t)|+(e^{M_{1}T}\|\rho_{0}\|_{L^{\infty}})^{2}+7M_{1}+4e^{M_{1}T}\|\rho_{0}\|_{L^{\infty}},$$
$$K_{4}=\max\limits_{t\in[0,T]}|h(t)|+\frac{1}{4}\left[2|a|+\int_{0}^{1}(u_{0,x}^{2}+\rho_{0}^{2})dx\right]^{2}
+(e^{M_{1}T}\|\rho_{0}\|_{L^{\infty}})^{2}.$$ By means of Gronwall's
inequality and the above inequality, we get
$$\|u(t,\cdot)\|_{H^{2}}^{2}+\|\rho(t,\cdot)\|_{H^{1}}^{2}\leq
e^{K_{3}t}(\|u_{0}\|_{H^{2}}^{2}+\|\rho_{0}\|_{H^{1}}^{2}+\frac{K_{4}}{K_{3}}).$$
The above inequality ensures that the solution $z$ does not blow-up
in finite time.

On the other hand, by Sobolev's imbedding theorem, we see that if
$$\liminf\limits_{t\rightarrow T}\inf\limits_{x\in
\mathbb{S}}{u_{x}(t,x)}=-\infty,$$ then the solution will blow up in
finite time. This completes the proof of the theorem.

\begin{remark3} Note that Theorem 3.2 and Theorem 3.5 show that
$$T(a,h(t),\|z_{0}\|_{H^{s}\times H^{s-1}})=T(a,h(t),\|z_{0}\|_{H^{s^{\prime}}
\times H^{s^{\prime}-1}})=T(a,h(t),\|z_{0}\|_{H^{2}\times H^{1}})$$
with $k=1$ for each \ $s,s^{\prime}>\frac{5}{2}.$ Furthermore, the
maximal existence time T of the family of solutions to (1.1) given
in Theorem 2.2 can be chosen independent of $s$. Moreover, Theorem
3.5 implies that
$$T(a,h(t),\|z_{0}\|_{H^{s}\times H^{s-1}})\leq T(a,h(t),\|z_{0}\|_{H^{2}\times
H^{1}})$$ with $k=\pm1$ for each $s\geq 2.$
\end{remark3}
\begin{remark3}
Note that Theorem 3.4 shows that $$T(a,h(t),\|z_{0}\|_{H^{s}\times
H^{s-1}})=T(a,h(t),\|z_{0}\|_{H^{s^{\prime}}\times
H^{s^{\prime}-1}})$$ with $k=-1$ for each \
$s,s^{\prime}>\frac{5}{2}.$ Moreover, the maximal existence time T
of the family of solutions to (1.1) given in Theorem 2.2 can be
chosen independent of $s$.
\end{remark3}
\section{Blow-up}
\newtheorem{theorem4}{Theorem}[section]
\newtheorem{lemma4}{Lemma}[section]
\newtheorem {remark4}{Remark}[section]
\newtheorem{corollary4}{Corollary}[section]

In this section, we discuss the blow-up phenomena of the system
(1.1) and prove that there exist strong solutions to (1.1) which do
not exist globally in time.

\begin{theorem4}Assume $k=1$. Let $z_{0}=\left(
                                                     \begin{array}{c}
                                                       u_{0} \\
                                                       \rho_{0} \\
                                                     \end{array}
                                                   \right)
\in H^s\times H^{s-1}, s\geq 2,$  and T be the maximal time of the
solution $z=\left(
                                    \begin{array}{c}
                                      u \\
                                      \rho\\
                                    \end{array}
                                  \right)$
to (1.1) with the initial data $z_0$. If $\rho_{0}\not\equiv0$ or
$u_{0}\not\equiv c$ for any $c\in\mathbb{R}$, and there exists a
point $x_{0}\in \mathbb{S},$ such that $\rho_{0}(x_{0})=0$, then the
corresponding solutions to (1.1) blow up in finite time.
\end{theorem4}
\textbf{Proof} We use the integrated representation (2.1). Let
$m(t)=u_{x}(t,q(t,x_{0}))$, $\gamma(t)=\rho (t,q(t,x_{0}))$, where
$q(t,x)$ is the solution of Eq.(3.9). By Eq.(3.9) we can obtain
$$\frac{dm}{dt}=(u_{tx}+uu_{xx})(t,q(t,x_{0})).$$
Evaluating (2.1) at $(t, q(t,x_{0}))$ we get
$$\frac{d}{dt}m(t)=\frac{1}{2}\gamma(t)^{2}-\frac{1}{2}m(t)^{2}+a.$$
Since $\gamma(0)=0$, we infer from Lemmas 3.5-3.6 that $\gamma(t)=0$
for all $t\in[0,T).$ Note that
$a=-\frac{1}{2}\int_{\mathbb{S}}(\rho_{0}^{2}+u_{0,x}^{2})dx<0$
since $\rho_{0}\not\equiv0$ or $u_{0}\not\equiv c$. Then we have
$\frac{d}{dt}m(t)\leq a < 0$. Thus, it follows that $m(t_{0})<0$ for
some $t_{0}\in(0,T).$ Solving the following inequality yields
$$\frac{d}{dt}m(t)\leq-\frac{1}{2}m(t)^{2}.$$ Therefore $$ 0 >
\frac{1}{m(t)}\geq \frac{1}{m(t_{0})}+\frac{1}{2}(t-t_{0}).$$ The
above inequality implies that $T< t_{0}-\frac{2}{m(t_{0})}$ and
$\lim\limits_{t\rightarrow T}m(t)= -\infty.$ In view of Theorem 3.5
and Remark 3.2, this completes the proof of the theorem.

\begin{corollary4}
Assume $k=1$. Let $z_{0}=\left(
                                                     \begin{array}{c}
                                                       u_{0} \\
                                                       \rho_{0} \\
                                                     \end{array}
                                                   \right)
\in H^s\times H^{s-1}, s\geq 2 ,$  and T be the maximal time of the
solution $z=\left(
                                    \begin{array}{c}
                                      u \\
                                      \rho \\
                                    \end{array}
                                  \right)
$ to (1.1) with the initial data $z_0$.  If $\rho_{0}$ is odd,
either $\rho_{0}\not\equiv 0$ or $u_{0}\not\equiv 0$ is odd, then
the corresponding solutions to (1.1) blow up in finite time.
\end{corollary4}
\textbf{Proof} Since $\rho_{0}$ is odd, $\rho_{0}(0)=0$.
$u_{0}\not\equiv 0$ being odd implies $u_{0}\not\equiv c$ for any
$c\in\mathbb{R}$. From Theorem 4.1 we can get the desired result.

\begin{theorem4}
Assume $k=-1$. Let $z_{0}=\left(
                                                     \begin{array}{c}
                                                       u_{0} \\
                                                       \rho_{0} \\
                                                     \end{array}
                                                   \right)
\in H^s\times H^{s-1}, s>\frac{5}{2},$  and T be the maximal time of
the solution $z=\left(
                                    \begin{array}{c}
                                      u \\
                                      \rho \\
                                    \end{array}
                                  \right)$
to (1.1) with the initial data $z_0$. The corresponding solutions to
(1.1) blow up in finite time if one of the following conditions
holds: (1) $a<0$, (2) $a> 0$ and there exists some $x_{0}\in
\mathbb{S}$ such that $u'_0(x_{0})<-\sqrt{2a},$ (3) $a=0$ and there
exists some $x_{0}\in \mathbb{S}$ such that
$u_{0}^{\prime}(x_{0})\leq 0$, $\rho_{0}(x_{0})\neq 0.$
\end{theorem4}
\textbf{Proof} Applying Remark 3.3 and a simply density argument, it
is clear that we may consider the case $s=3.$ Define now
$$ m(t):=\min_{x\in \mathbb{S}}\{u_{x}(t,x)\}, \ \ t\in [0,T)$$ and
let $\xi(t)\in \mathbb{S}$ be a point where this minimum is attained
by Lemma 3.3. It follows that
$$m(t)=u_{x}(t,\xi(t)).$$ Clearly $u_{xx}(t,\xi(t))=0$ since $u(t,\cdot)\in
H^{3}(\mathbb{S})\subset C^{2}(\mathbb{S}).$ Using the integrated
representation (2.1) and evaluating it at $(t,\xi(t))$, we obtain
$$\frac{d}{dt}m(t)\leq-\frac{1}{2}m(t)^{2}+a.$$

Let (1) hold. Note that $\frac{dm(t)}{dt}\leq a$. It then follows
that there is a point $x_{0}\in \mathbb{S}$ such that $m(t_{0})< 0$.
Solving the following inequality
$$\frac{d}{dt}m(t)\leq-\frac{1}{2}m(t)^{2},$$ we obtain $$ 0 >
\frac{1}{m(t)}\geq \frac{1}{m(t_{0})}+\frac{1}{2}(t-t_{0}).$$ This
implies that $T< t_{0}-\frac{2}{m(t_{0})}$ and
$\lim\limits_{t\rightarrow T}m(t)= -\infty.$

Let (2) hold. Note that if $m(0)=u'_0(\xi(0))\leq
u'_0(x_{0})<-\sqrt{2a}$, then $m(t)<-\sqrt{2a}$ for all $t\in[0,T)$.
From the above inequality we obtain
\begin{eqnarray*}
\frac{m(0)+\sqrt{2a}}{m(0)-\sqrt{2}K}e^{\sqrt{2a}\
t}-1\leq\frac{2\sqrt{2a}}{m(t)-\sqrt{2a}}\leq 0.\end{eqnarray*}
Since $0<\frac{m(0)+\sqrt{2a}}{m(0)-\sqrt{2a}}<1,$ there exists
$$ 0<T\leq\frac{1}{\sqrt{2a}}\ln\frac{m(0)-\sqrt{2a}}{m(0)+\sqrt{2}a},$$
 such that $\lim_{t\rightarrow T}m(t)=-\infty.$
Theorem 3.5 and Remark 3.2 imply that the corresponding solution to
(2.3) blows up in finite time if condition (1) or condition (2)
holds.

Let (3) hold. We use the integrated representation (2.1). Let
$h(t)=u_{x}(t,q(t,x_{0}))$, $\gamma(t)=\rho (t,q(t,x_{0}))$, where
$q(t,x)$ is the solution of Eq.(3.9). By Eq.(3.9) we can obtain
$$\frac{dh}{dt}=(u_{tx}+uu_{xx})(t,q(t,x_{0})).$$
Evaluating (2.1) at $(t, q(t,x_{0}))$ we get
$$\frac{d}{dt}h(t)=-\frac{1}{2}\gamma(t)^{2}-\frac{1}{2}h(t)^{2}.$$
By $\gamma(0)\neq0$, we infer from Lemmas 3.5-3.6 that
$\gamma(t)\neq0$ for all $t\in[0,T).$ Since $h(0)\leq 0$ and
$\frac{d}{dt}h(t)<0$, it follows that $h(t_{0})< 0$ for some
$t_{0}\in [0,T).$ Solving the following inequality
$$\frac{d}{dt}h(t)\leq-\frac{1}{2}h(t)^{2},$$ we obtain $$ 0 >
\frac{1}{h(t)}\geq \frac{1}{h(t_{0})}+\frac{1}{2}(t-t_{0}).$$ This
implies that $T< t_{0}-\frac{2}{h(t_{0})}$ and
$\lim\limits_{t\rightarrow T}h(t)= -\infty.$ In view of Theorem 3.5
and Remark 3.2, this completes the proof of the theorem.

\begin{corollary4}
Assume $k=-1$. Let $z_{0}=\left(
                                                     \begin{array}{c}
                                                       u_{0} \\
                                                       \rho_{0} \\
                                                     \end{array}
                                                   \right)
\in H^s\times H^{s-1}, s>\frac{5}{2},$  and T be the maximal time of
the solution $z=\left(
                                    \begin{array}{c}
                                      u \\
                                      \rho \\
                                    \end{array}
                                  \right)$
to (2.3) with the initial data $z_0$. The corresponding solution to
(2.3) blows up in finite time if one of the following conditions
holds: (1) $a> 0$ and $u'_0(0)<-\sqrt{2a},$ (2) $a=0$ and
$u_{0}^{\prime}(0)\leq 0$, $\rho_{0}(0)\neq 0.$
\end{corollary4}

\section{Global Existence}
\newtheorem{theorem5}{Theorem}[section]
\newtheorem{lemma5}{Lemma}[section]
\newtheorem {remark5}{Remark}[section]
\newtheorem{corollary5}{Corollary}[section]
In this section, we will present a global existence result, which
improves considerably the recent results in \cite{C-I, MW}.
\begin{theorem5}Assume $k=1$. Let $z_{0}=\left(
                                                     \begin{array}{c}
                                                       u_{0} \\
                                                       \rho_{0} \\
                                                     \end{array}
                                                   \right)
\in H^{s}\times H^{s-1},$ where $s=2$ or $s\geq 3$ and T be the
maximal time of the solution $z=\left(
                                    \begin{array}{c}
                                      u \\
                                      \rho \\
                                    \end{array}
                                  \right)
$ to (1.1) with the initial data $z_0$. If $\rho_{0}(x)\neq 0$ for
all $x\in\mathbb{S}$, then the corresponding solutions $z$ exist
globally in time.
\end{theorem5}
\textbf{Proof}  By Lemma 3.5, we know that $q(t,\cdot)$ is an
increasing diffeomorphism of $\mathbb{R}$ with
$$
q_{x}(t,x)=exp\left(\int_{0}^{t}u_{x}(s,q(s,x))ds\right)>0, \ \
\forall \ (t,x)\in [0,T)\times \mathbb{R}.$$ Moreover,
\begin{eqnarray}
\inf_{y\in\mathbb{S}}u_{x}(t,y)=\inf_{x\in\mathbb{R}}u_{x}(t,q(t,x)),
 \ \forall \ t\in[0,T).
\end{eqnarray}
Set $M(t,x)=u_{x}(t,q(t,x))$ and $\alpha(t,x)=\rho(t,q(t,x))$ for
$t\in [0,T)$ and $x\in\mathbb{R}$. By (1.1) and Eq.(3.9), we have
\begin{eqnarray}
\frac{\partial M}{\partial t}=(u_{tx}+uu_{xx})(t,q(t,x))\ \
\text{and} \ \ \frac{\partial \alpha}{\partial t}=-\alpha M.
\end{eqnarray}
Evaluating (2.1) at $(t,q(t,x))$ we get
\begin{equation}
\partial_{t}M(t,x)=-\frac{1}{2}M(t,x)^{2}+\frac{1}{2}\alpha(t,x)^{2}+a.
\end{equation}

By Lemmas 3.5-3.6, we know that $\alpha(t,x)$ has the same sign with
$\alpha(0,x)=\rho_{0}(x)$ for every $x\in \mathbb{R}$. Moreover,
there is a constant $\beta>0$ such that
$\inf\limits_{x\in\mathbb{R}}|\alpha(0,x)|=\inf\limits_{x\in\mathbb{S}}|\rho_{0}(x)|\geq\beta>0$
since $\rho_{0}(x)\neq 0$ for all $x\in\mathbb{S}$ and $\mathbb{S}$
is a compact set. Thus,
$$\alpha(t,x)\alpha(0,x)>0, \ \ \ \forall x\in \mathbb{R}.$$Next, we
consider the following Lyapunov function first introduced in
\cite{C-I}.
\begin{equation}
w(t,x)=\alpha(t,x)\alpha(0,x)+\frac{\alpha(0,x)}{\alpha(t,x)}(1+M^{2}),
\ \ \ (t,x)\in[0,T)\times\mathbb{R}.
\end{equation}
By Sobolev's imbedding theorem, we have
\begin{align}
0&<w(0,x)=\alpha(0,x)^{2}+1+M(0,x)^{2}\\
&\nonumber=\rho_{0}(x)^{2}+1+u_{0,x}(x)^{2}\\
&\nonumber\leq1+\max\limits_{x\in\mathbb{S}}(\rho_{0}(x)^{2}+u_{0,x}(x)^{2}):=C_{1}.
\end{align}
 Differentiating
(5.4) with respect to $t$ and using (5.2)-(5.3), we obtain

\begin{align}
\nonumber \frac{\partial w}{\partial
t}(t,x)&=\frac{\alpha(0,x)}{\alpha(t,x)}M(t,x)\left(2a+1\right)\\
&\nonumber
\leq|1+2a|\frac{\alpha(0,x)}{\alpha(t,x)}(1+M^{2})\\
&\nonumber\leq|1+2a|w(t,x).
\end{align}
By Gronwall's inequality, the above inequality and (5.5), we have
$$
w(t,x)\leq w(0,x)e^{|1+2a|t}\leq C_{1}e^{|1+2a|t}
$$
for all $(t,x)\in [0,T)\times\mathbb{R}.$ On the other
hand,$$w(t,x)\geq2\sqrt{\alpha^{2}(0,x)(1+M^{2})}\geq2\beta|M(t,x)|,
\ \ \ \forall \ \ (t,x)\in [0,T)\times\mathbb{R}.$$ Thus,
$$
M(t,x)\geq-\frac{1}{2\beta}w(t,x)\geq-\frac{1}{2\beta}C_{1}e^{|1+2a|t}
$$
for all $(t,x)\in [0,T)\times\mathbb{R}$. Then by (5.1) and the
above inequality, we have
$$\lim\limits_{t\rightarrow
T}\inf\limits_{y\in\mathbb{S}}u_{x}(t,y)=\lim\limits_{t\rightarrow
T}\inf\limits_{x\in\mathbb{R}}u_{x}(t,q(t,x))
\geq-\frac{1}{2\beta}C_{1}e^{|1+2a|t}.
$$
This completes the proof by using Theorem 3.2 and Remark 3.2.

\bigskip

\noindent\textbf{Acknowledgments} This work was partially supported
by NNSFC (No. 10971235), RFDP (No. 200805580014), NCET-08-0579 and
the key project of Sun Yat-sen University.

\end{document}